\documentstyle[12pt]{article}

\font\twlgot =eufm10 scaled \magstep1
\font\egtgot =eufm8
\font\sevgot =eufm7
\font\twlmsb =msbm10 scaled \magstep1
\font\egtmsb =msbm8
\font\sevmsb =msbm7

\newfam\gotfam

\textfont\gotfam\twlgot
\scriptfont\gotfam\egtgot
\scriptscriptfont\gotfam\sevgot

\newfam\msbfam
\textfont\msbfam\twlmsb
\scriptfont\msbfam\egtmsb
\scriptscriptfont\msbfam\sevmsb
\def\Bbb{\protect\pBbb}
\def\pBbb{\relax\ifmmode\expandafter\Bb\else\typeout{You cann't use
Bbb in text mode}\fi}
\def\Bb #1{{\fam\msbfam\relax#1}}

\def\thebibliography#1{\section*{References}\list
  {[\arabic{enumi}]}{\settowidth\labelwidth{#1}\leftmargin\labelwidth
    \advance\leftmargin\labelsep
    \usecounter{enumi}}
    \def\newblock{\hskip .11em plus .33em minus .07em}
    \sloppy\clubpenalty4000\widowpenalty4000
    \sfcode`\.=1000\relax}

\newcommand{\beq}{\begin{equation}}
\newcommand{\eeq}{\end{equation}}
\newcommand{\ben}{\begin{eqnarray}}
\newcommand{\een}{\end{eqnarray}}
\newcommand{\be}{\begin{eqnarray*}}
\newcommand{\ee}{\end{eqnarray*}}
\newcommand{\bea}{\begin{eqalph}}
\newcommand{\eea}{\end{eqalph}}

\newcommand{\al}{\alpha}

\newcommand{\bt}{\beta}
\newcommand{\dl}{\delta}
\newcommand{\la}{\lambda}

\newcommand{\f}{\phi}

\newcommand{\Om}{\Omega}
\newcommand{\m}{\mu}

\newcommand{\vt}{\vartheta}
\newcommand{\vf}{\varphi}

\newcommand{\w}{\wedge}

\newcommand{\dr}{\partial}

\newcounter{eqalph}
\newcounter{equationa}
\newcounter{theorem}
\newcounter{remark}
\newcounter{proposition}
\newcounter{lemma}
\newcounter{corollary}
\newcounter{definition}

\newenvironment{eqalph}{\stepcounter{equation}
\setcounter{equationa}{\value{equation}}
\setcounter{equation}{0}

\begin{eqnarray}}{\end{eqnarray}\setcounter{equation}{\value{equationa}}}

\def\theremark{\arabic{remark}}

\def\thedefinition{\arabic{definition}}

\newenvironment{proof}{\noindent 
{\it Proof.}}{\hfill{\footnotesize\it QED} \bigskip}

\newenvironment{theo}{\refstepcounter{definition} 
\bigskip\noindent{\it Theorem \thedefinition.}}{\medskip}

\newcommand{\mar}[1]{}

\begin{document}
\hbox{}

{\parindent=0pt

{\large\bf Action-angle coordinates around a noncompact 
invariant manifold of a completely integrable Hamiltonian system}
\bigskip

{\sc E.Fiorani}\footnote{E-mail: fiorani@mat.unimi.it}

{\sl Dipartimento di Matematica "F.Enriques", Universit\'a di Milano, 
20133, Milano, Italy}

{\sc G.Giachetta}\footnote{E-mail: giovanni.giachetta@unicam.it},

{\sl Dipartimento di Matematica e Fisica, Universit\'a di Camerino,
62032 Camerino (MC), Italy}

\medskip

{\sc G. Sardanashvily}\footnote{E-mail:
sard@grav.phys.msu.su; URL: http://webcenter.ru/$\sim$sardan/}

{\sl Department of Theoretical Physics,
Moscow State University, 117234 Moscow, Russia}

\bigskip
\bigskip
{\small

{\bf Abstract.}
A trivial bundle of regular connected invariant
manifolds of a completely integrable
Hamiltonian system can be provided with
action-angle coordinates.   
\bigskip


}}

\bigskip\bigskip

By virtue of the 
classical theorem \cite{arn,laz}, an autonomous completely
integrable Hamiltonian system admits
action-angle coordinates around a connected regular compact invariant manifold.
We show that 
there is a system of action-angle coordinates on an open neighbourhood
$U$ of a connected regular
invariant manifold $M$ 
if Hamiltonian vector fields of first integrals on $U$ are complete
and the foliation of $U$ by invariant manifolds is trivial.
If $M$ is a compact regular invariant manifold, these conditions always hold.

\begin{theo} \label{z8} \mar{z8}
Let $M$ be a connected invariant manifold of a completely integrable
Hamiltonian system $\{F_\la\}$, $\la=1,\ldots,n$, on a symplectic manifold
$(Z,\Om)$. Let $U$ be an open neighbourhood of $M$ such that: (i) $\{F_\la\}$ have
no critical points in $U$, (ii) the Hamiltonian vector fields of 
the first integrals $F_\la$ on $U$ are complete,
and (iii) the submersion $\times F_\la: U\to \Bbb R^n$
is a trivial bundle of invariant manifolds 
over a domain $V\subset \Bbb R^n$. 
Then $U$ is isomorphic
to the symplectic annulus 
\mar{z10}\beq
W=V\times(\Bbb R^{n-m}\times T^m), \label{z10}
\eeq
provided with the action-angle coordinates 
\be
(I_1,\ldots,I_n; x^1,\ldots, x^{n-m}; \f^1,\ldots,\f^m) 
\ee
such that the symplectic form on $W$ reads
\be
\Om=dI_a\w dx^a +dI_i\w d\f^i,
\ee
and the first integrals $F_\la$ depend only on  
the action coordinates $(I_\al)$.
\end{theo}

\begin{proof}
In accordance with the well-known theorem \cite{arn}, 
the invariant
manifold $M$ is diffeomorphic to the product $\Bbb R^{n-m}\times T^m$,
which is the group space of the quotient $G=\Bbb R^n/\Bbb Z^m$ 
of the group $\Bbb R^n$ generated by Hamiltonian vector fields
$\vt_\la$ of first integrals $F_\la$ on $M$. 
Namely, $M$
is provided with the group space coordinates $(y^\la)=(s^a,\vf^i)$ 
where $\vf^i$ are linear functions 
of parameters $s^\la$ along integral curves of the 
Hamiltonian vector fields $\vt_\la$ on $U$. Let
$(J_\la)$ be coordinates on
$V$ which 
are values of first integrals $F_\la$. Let us choose a trivialization
of the fiber bundle $U\to V$ seen as a principal bundle with the
structure group $G$. We fix its global section $\chi$.
Since parameters $s^\la$ are given up to a shift, let us
provide each fiber $M_J$, $J\in V$, with the group space
coordinates $(y^\la)$ 
centered at the point 
$\chi(J)$. Then 
$(J_\la;y^\la)$ are bundle coordinates
on the annulus $W$ (\ref{z10}). 
Since $M_J$ are Lagrangian manifolds, 
the symplectic form $\Om$ on $W$
is given relative to the
bundle coordinates $(J_\la;y^\la)$ by the expression
\mar{ac1}\beq
\Om=\Om^{\al\bt}dJ_\al\w dJ_\bt + \Om^\al_\bt dJ_\al\w dy^\bt. \label{ac1}
\eeq
 By the very definition of coordinates $(y^\la)$, the
Hamiltonian vector fields $\vt_\la$ of first integrals take the
coordinate form $\vt_\la=\vt_\la^\al(J_\m)\dr_\al$. Moreover, since
the cyclic group $S^1$ can not act transitively on $\Bbb R$, we have
\mar{ww25}\beq
\vt_a=\dr_a +\vt_a^i(J_\la)\dr_i, \qquad \vt_i=\vt_i^k(J_\la)\dr_k. \label{ww25}
\eeq
The Hamiltonian vector fields $\vt_\la$ obey the relations
\mar{ww22}\beq
\vt_\la\rfloor\Om=-dJ_\la,\qquad
\Om^\al_\bt \vt^\bt_\la=\dl^\al_\la. \label{ww22}
\eeq
It follows that $\Om^\al_\bt$ is a nondegenerate matrix and 
$\vt^\al_\la=(\Om^{-1})^\al_\la$, i.e., the functions $\Om^\al_\bt$
depend only on coordinates $J_\la$.  
A substitution of (\ref{ww25}) into (\ref{ww22}) results in the equalities 
\mar{ww30,1}\ben
&& \Om^a_b=\dl^a_b, \qquad \vt_a^\la\Om^i_\la=0, \label{ww30}\\
&& \vt^k_i\Om^j_k=\dl^j_i, \qquad \vt^k_i\Om^a_k=0. \label{ww31}
\een
The first of the equalities (\ref{ww31}) shows that the matrix $\Om^j_k$
is nondegenerate, and so is the matrix $\vt^k_i$. Then the second 
one gives $\Om^a_k=0$.

By virtue of
the well-known K\"unneth formula for the de Rham cohomology of a product of manifolds,
the closed form $\Om$ (\ref{ac1}) on $W$ (\ref{z10})
is exact, i.e., $\Om=d\Xi$ where $\Xi$  reads
\be
\Xi=\Xi^\al(J_\la,y^\la)dJ_\al + \Xi_i(J_\la) d\vf^i + 
\dr_\al\Phi(J_\la,y^\la)dy^\al, 
\ee
where $\Phi$ is a function on $W$.
Taken up to an exact form, $\Xi$ is brought into the form
\mar{ac2}\beq
\Xi=\Xi'^\al(J_\la,y^\la)dJ_\al + \Xi_i(J_\la) d\vf^i. \label{ac2}
\eeq
Owing to the fact that
components of $d\Xi=\Om$ are independent of $y^\la$ and obey the equalities
(\ref{ww30}) -- (\ref{ww31}), we obtain the following.
 
(i) $\Om^a_i=-\dr_i\Xi'^a +\dr^a\Xi_i=0$. It follows that $\dr_i\Xi'^a$ is 
independent of $\vf^i$, i.e., $\Xi'^a$ is affine in $\vf^i$ and, consequently,
is independent of $\vf^i$ since $\vf^i$ are cyclic coordinate. Hence, 
$\dr^a\Xi_i=0$, i.e., $\Xi_i$ is a function only of coordinates $J_j$.

(ii) $\Om^k_i=-\dr_i\Xi'^k +\dr^k\Xi_i$. Similarly to item (i), one shows that $\Xi'^k$
is independent of $\vf^i$ and $\Om^k_i=\dr^k\Xi_i$, i.e., 
$\dr^k\Xi_i$ is a nondegenerate matrix.

(iii) $\Om^a_b=-\dr_b\Xi'^a=\dl^a_b$. Hence, $\Xi'^a=-s^a+D^a(J_\la)$.

(iv) $\Om^i_b=-\dr_b\Xi'^i$, i.e., $\Xi'^i$ is affine in $s^a$.

In view of items (i) -- (iv), the Liouville form $\Xi$ (\ref{ac2}) reads
\be
\Xi=x^adJ_a + [D^i(J_\la) + B^i_a(J_\la)s^a]dJ_i + \Xi_i(J_j) d\vf^i,
\ee
where we put
\mar{ee1}\beq
x^a=-\Xi'^a=s^a-D^a(J_\la). \label{ee1}
\eeq
Since the matrix $\dr^k\Xi_i$ is nondegenerate,
one can introduce new coordinates $I_i=\Xi_i(J_j)$, $I_a=J_a$. Then we have
\be
\Xi=-x^adI_a + [D'^i(I_\la) + B'^i_a(I_\la)s^a]dI_i + I_i d\vf^i.
\ee
Finally, put 
\mar{ee2}\beq
\f^i=\vf^i-[D'^i(I_\la) + B'^i_a(I_\la)s^a] \label{ee2}
\eeq
in order to obtain the desired action-angle coordinates
\be
I_a=J_a, \qquad I_i(J_j), \qquad x^a(J_\la,s^a), \qquad \f^i(J_\la,y^\la).
\ee
These are bundle coordinates on $U\to V$ where the 
coordinate shifts (\ref{ee1}) -- (\ref{ee2}) correspond to a choice of another
trivialization of $U\to V$.
\end{proof}

In particular, Theorem \ref{z8} enables one to introduce action-angle coordinates
for time-dependent completely integrable Hamiltonian 
systems whose invariant manifolds are
never compact because of the time axis \cite{epr}. Another application of Theorem
\ref{z8} is Jacobi fields of a completely integrable system. They also make up
a completely integrable system whose invariant manifolds are never compact \cite{epr1}.

\end{document}